\input amstex
\documentstyle{amsppt}
%
\catcode`@=11
\redefine\output@{%
  \def\break{\penalty-\@M}\let\par\endgraf
  \ifodd\pageno\global\hoffset=105pt\else\global\hoffset=8pt\fi  
  \shipout\vbox{%
    \ifplain@
      \let\makeheadline\relax \let\makefootline\relax
    \else
      \iffirstpage@ \global\firstpage@false
        \let\rightheadline\frheadline
        \let\leftheadline\flheadline
      \else
        \ifrunheads@ 
        \else \let\makeheadline\relax
        \fi
      \fi
    \fi
    \makeheadline \pagebody \makefootline}%
  \advancepageno \ifnum\outputpenalty>-\@MM\else\dosupereject\fi
}
\def\Beta{\mathchar"0\hexnumber@\rmfam 42}
\catcode`\@=\active
\nopagenumbers
\chardef\textvolna='176
\def\negskp{\hskip -2pt}
\def\Ker{\operatorname{Ker}}

\def\Sym{\operatorname{Sym}}
\chardef\degree="5E
\def\blue#1{#1}

\catcode`#=11\def\diez{#}\catcode`#=6
\catcode`&=11\catcode`&=4
\catcode`_=11\def\podcherkivanie{_}\catcode`_=8
\catcode`~=11\def\volna{~}\catcode`~=\active
\def\mycite#1{\cite{\blue{#1}}\immediate\special{ps:
     ShrHPSdict begin /ShrBORDERthickness 0 def}}
\def\myciterange#1#2#3#4{\cite{\blue{#2#3#4}}\immediate\special{ps:
     ShrHPSdict begin /ShrBORDERthickness 0 def}}
\def\mytag#1{%
    \tag#1}
\def\mythetag#1{\thetag{\blue{#1}}\immediate\special{ps:
     ShrHPSdict begin /ShrBORDERthickness 0 def}}
\def\myrefno#1{\no#1}
\def\myhref#1#2{\blue{#2}\immediate\special{ps:
     ShrHPSdict begin /ShrBORDERthickness 0 def}}
\def\myEarXivlink{\myhref{http://arXiv.org}{http:/\negskp/arXiv.org}}

\def\mytheorem#1{\csname proclaim\endcsname{Theorem #1}}
\def\mytheoremwithtitle#1#2{\csname proclaim\endcsname{Theorem #1#2}}
\def\mythetheorem#1{\blue{#1}\immediate\special{ps:
     ShrHPSdict begin /ShrBORDERthickness 0 def}}
\def\mylemma#1{\csname proclaim\endcsname{Lemma #1}}
\def\mylemmawithtitle#1#2{\csname proclaim\endcsname{Lemma #1#2}}

\def\mycorollary#1{\csname proclaim\endcsname{Corollary #1}}

\def\mydefinition#1{\definition{Definition #1}}

\def\myconjecture#1{\csname proclaim\endcsname{Conjecture #1}}
\def\myconjecturewithtitle#1#2{\csname proclaim\endcsname{Conjecture #1#2}}

\def\myproblem#1{\csname proclaim\endcsname{Problem #1}}
\def\myproblemwithtitle#1#2{\csname proclaim\endcsname{Problem #1#2}}
\def\mytheproblem#1{\blue{#1}\immediate\special{ps:
     ShrHPSdict begin /ShrBORDERthickness 0 def}}

\font\eightcyr=wncyr8
\pagewidth{360pt}
\pageheight{606pt}
\topmatter
\title
Inverse problems associated with perfect cuboids.
\endtitle
\rightheadtext{Inverse problems associated with perfect cuboids.}
\author
John Ramsden, Ruslan Sharipov
\endauthor
\address CSR Plc, Cambridge Business Park, Cambridge, CB4 0WZ, UK
\endaddress
\email\myhref{mailto:jhnrmsdn\@yahoo.co.uk}{jhnrmsdn\@yahoo.co.uk}
\endemail
\address Bashkir State University, 32 Zaki Validi street, 450074 Ufa, Russia
\endaddress
\email\myhref{mailto:r-sharipov\@mail.ru}{r-sharipov\@mail.ru}
\endemail
\abstract
    A perfect cuboid is a rectangular parallelepiped with integer edges, 
integer face diagonals, and integer space diagonal. Such cuboids have not yet 
been found, but nor has their existence been disproved. Perfect cuboids are 
described by a certain system of Diophantine equations possessing an intrinsic 
$S_3$ symmetry. Recently these equations were factorized with respect to this 
$S_3$ symmetry and the factor equations were transformed into $E$-form. As 
appears, the transformed factor equations are explicitly solvable. Based on this 
solution, polynomial inverse problems are formulated in the present paper.
\endabstract
\endtopmatter
\TagsOnRight
\document

%
%
\head
1. Introduction.
\endhead
     Perfect cuboids are described by the following four polynomial equations:
$$
\xalignat 2
&\hskip -2em
x_1^2+x_2^2+x_3^2-L^2=0,
&&x_2^2+x_3^2-d_1^{\kern 1pt 2}=0,\\
\vspace{-1.7ex}
\mytag{1.1}\\
\vspace{1ex}
\vspace{-1.7ex}
&\hskip -2em
x_3^2+x_1^2-d_2^{\kern 1pt 2}=0,
&&x_1^2+x_2^2-d_3^{\kern 1pt 2}=0.
\endxalignat
$$
Here $x_1$, $x_2$, $x_3$ are edges of a cuboid, $d_1$, $d_2$, $d_3$ are its face 
diagonals, and $L$ is its space diagonal. For the history of perfect cuboids the
reader is referred to \myciterange{1}{1}{--}{44}.\par
     Recently in \mycite{45} the symmetry approach to the equations \mythetag{1.1}
was initiated. It is based on an intrinsic $S_3$ symmetry of these equations. 
Indeed, let's consider the following action of the group $S_3$ upon the variables 
$x_1$, $x_2$, $x_3$, $d_1$, $d_2$, $d_3$ and $L$:
$$
\xalignat 3
&\hskip -2em
\sigma(x_i)=x_{\sigma i},
&&\sigma(d_i)=d_{\sigma i},
&&\sigma(L)=L.
\mytag{1.2}
\endxalignat
$$ 
The first equation \mythetag{1.1} is invariant with respect to the transformations
\mythetag{1.2}. The other three equations are not invariant, but the system as a whole
is again invariant, i\.\,e\. it possesses $S_3$ symmetry based on the transformations 
\mythetag{1.2}. In \mycite{46} the equations \mythetag{1.1} were factorized with 
respect to their $S_3$ symmetry and the following system of eight factor equations 
was derived:
$$
\align
\hskip -2em
x_1^2+x_2^2+x_3^2-L^2=0,
\mytag{1.3}\\
\vspace{1ex}
\hskip -2em
(x_2^2+x_3^2-d_1^{\kern 1pt 2})
+(x_3^2+x_1^2-d_2^{\kern 1pt 2})+(x_1^2+x_2^2-d_3^{\kern 1pt 2})=0,
\mytag{1.4}\\
\vspace{1ex}
\hskip -2em
d_1\,(x_2^2+x_3^2-d_1^{\kern 1pt 2})
+d_2\,(x_3^2+x_1^2-d_2^{\kern 1pt 2})+d_3\,(x_1^2+x_2^2-d_3^{\kern 1pt 2})=0,
\mytag{1.5}\\
\endalign
$$
$$
\align
&\hskip -2em
x_1\,(x_2^2+x_3^2-d_1^{\kern 1pt 2})
+x_2\,(x_3^2+x_1^2-d_2^{\kern 1pt 2})+x_3\,(x_1^2+x_2^2-d_3^{\kern 1pt 2})=0,
\mytag{1.6}\\
\vspace{2ex}
&\hskip -2em
\aligned
&x_1\,d_1\,(x_2^2+x_3^2-d_1^{\kern 1pt 2})
+x_2\,d_2\,(x_3^2+x_1^2-d_2^{\kern 1pt 2})\,+\\
&\hskip 14.6em +\,x_3\,d_3\,(x_1^2+x_2^2-d_3^{\kern 1pt 2})=0,
\endaligned
\mytag{1.7}\\
\vspace{2ex}
&\hskip -2em
x_1^2\,(x_2^2+x_3^2-d_1^{\kern 1pt 2})
+x_2^2\,(x_3^2+x_1^2-d_2^{\kern 1pt 2})+x_3^2\,(x_1^2+x_2^2-d_3^{\kern 1pt 2})=0,
\mytag{1.8}\\
\vspace{2ex}
&\hskip -2em
d_1^{\kern 1pt 2}\,(x_2^2+x_3^2-d_1^{\kern 1pt 2})
+d_2^{\kern 1pt 2}\,(x_3^2+x_1^2-d_2^{\kern 1pt 2})+d_3^{\kern 1pt 2}
\,(x_1^2+x_2^2-d_3^{\kern 1pt 2})=0,
\mytag{1.9}\\
\vspace{2ex}
&\hskip -2em
\aligned
&x_1^2\,d_1^{\kern 1pt 2}\,(x_2^2+x_3^2
-d_1^{\kern 1pt 2})+x_2^2\,d_2^{\kern 1pt 2}\,(x_3^2+x_1^2
-d_2^{\kern 1pt 2})\,+\\
&\hskip 14.6em +\,x_3^2\,d_3^{\kern 1pt 2}\,(x_1^2
+x_2^2-d_3^{\kern 1pt 2})=0.
\endaligned
\mytag{1.10}
\endalign
$$
Each solution of the equations \mythetag{1.1} is a solution of the equations
\mythetag{1.3}, \mythetag{1.4}, \mythetag{1.5}, \mythetag{1.6}, 
\mythetag{1.7}, \mythetag{1.8}, \mythetag{1.9}, \mythetag{1.10}. But generally
speaking, the converse is not true. Fortunately, in \mycite{47} the following
theorem was proved. 
\mytheorem{1.1} Each integer or rational solution of the equations \mythetag{1.3} 
through \mythetag{1.10} such that $x_1>0$, $x_2>0$, $x_3>0$, $d_1>0$, $d_2>0$, 
and $d_3>0$ is an integer or rational solution for the equations \mythetag{1.1}.
\endproclaim
     Due to Theorem~\mythetheorem{1.1} the factor equations \mythetag{1.3} 
through \mythetag{1.10} can be applied for studying perfect cuboids.\par
     Note that each of the equations \mythetag{1.3} through \mythetag{1.10} is
invariant with respect to the transformations \mythetag{1.2}. The left hand sides
of these equations are multisymmetric polynomials in the sense of the following
definition. 
\mydefinition{1.1} A polynomial $p\in\Bbb Q[x_1,x_2,x_3,d_1,d_2,d_3,L]$ is called 
multisymmetric if it is invariant with respect to the action \mythetag{1.2} of the 
group $S_3$. 
\enddefinition
     Multisymmetric polynomials constitute a subring within the polynomial ring
$\Bbb Q[x_1,x_2,x_3,d_1,d_2,d_3,L]$. For the sake of brevity we introduce the 
matrix
$$
M=\Vmatrix x_1 & x_2 &x_3\\
\vspace{1ex}
d_1 & d_2 & d_3\endVmatrix
$$
and denote $Q[x_1,x_2,x_3,d_1,d_2,d_3,L]=\Bbb Q[M,L]$. Similarly, the subring
of multisymmetric polynomials is denoted through $\Sym\!\Bbb Q[M,L]$. For the 
general theory of multisymmetric polynomials the reader is referred to 
\myciterange{48}{48}{--}{68}. For our purposes we need the following theorem 
from this theory.
\mytheorem{1.2} Each multisymmetric polynomial $p\in\Sym\!\Bbb Q[x_1,x_2,x_3,d_1,
d_2,d_3,L]$ can be expressed through the following elementary multisymmetric 
polynomials:
$$
\gather
\hskip -2em
\aligned
&e_{\sssize [1,0]}=x_1+x_2+x_3,\\
&e_{\sssize [2,0]}=x_1\,x_2+x_2\,x_3+x_3\,x_1,\\
&e_{\sssize [3,0]}=x_1\,x_2\,x_3,
\endaligned
\kern 3em
\aligned
&e_{\sssize [0,1]}=d_1+d_2+d_3,\\
&e_{\sssize [0,2]}=d_1\,d_2+d_2\,d_3+d_3\,d_1,\\
&e_{\sssize [0,3]}=d_1\,d_2\,d_3,
\endaligned
\quad
\mytag{1.11}\\
\vspace{2ex}
\hskip -2em
\aligned
&e_{\sssize [2,1]}=x_1\,x_2\,d_3+x_2\,x_3\,d_1+x_3\,x_1\,d_2,\\
&e_{\sssize [1,1]}=x_1\,d_2+d_1\,x_2+x_2\,d_3+d_2\,x_3+x_3\,d_1+d_3\,x_1,\\
&e_{\sssize [1,2]}=x_1\,d_2\,d_3+x_2\,d_3\,d_1+x_3\,d_1\,d_2.
\endaligned
\mytag{1.12}
\endgather
$$
\endproclaim
     Theorem~\mythetheorem{1.2} is known as the Fundamental Theorem on
Elementary Multisymmetric Polynomials. Its proof can be found in \mycite{54}.
If we denote 
$$
\Bbb Q[E_{10},E_{20},E_{30},E_{01},E_{02},E_{03},E_{21},E_{11},E_{12},L]
=\Bbb Q[E,L],
$$
then Theorem~\mythetheorem{1.2} means that for each multisymmetric 
polynomial $p$ from the ring $\Sym\!\Bbb Q[M,L]$ there is some polynomial 
$q\in\Bbb Q[E,L]$ of ten independent variables such that $p$ is produced 
from $q$ by substituting the elementary multisymmetric polynomials 
\mythetag{1.11} and \mythetag{1.12} for $E_{10}$, $E_{20}$, $E_{30}$, $E_{01}$, 
$E_{02}$, $E_{03}$, $E_{21}$, $E_{11}$, and $E_{12}$ into its arguments. 
The substitution procedure can be understood as a mapping:
$$
\hskip -2em
\varphi\!:\,\Bbb Q[E,L]\to\Sym\!\Bbb Q[M,L].
\mytag{1.13}
$$
The mapping \mythetag{1.13} is a ring homomorphism. It is surjective, which follows 
from Theorem~\mythetheorem{1.2}, but it is not bijective, i\.\,e\. it has a 
nonzero kernel $K=\Ker\varphi$. The kernel $K$ is an ideal in the ring $\Bbb Q[E,L]$.
It was calculated as
$$
\hskip -2em
K=\bigl<q_{\kern 1pt 1},\,\ldots,\,q_{\kern 1pt 7}\bigr>,
\mytag{1.14}
$$
where $q_{\kern 1pt 1},\,\ldots,\,q_{\kern 1pt 7}$ are seven polynomials being 
a basis of $K$. These polynomials are given by the explicit formulas \thetag{2.4} 
through \thetag{2.10} in \mycite{46}. The ideal \mythetag{1.14} has a Gr\"obner 
basis comprising fourteen polynomials:
$$
\hskip -2em
K=\bigl<\tilde q_{\kern 1pt 1},\tilde q_{\kern 1pt 2},\tilde q_{\kern 1pt 3},
\tilde q_{\kern 1pt 4},\tilde q_{\kern 1pt 5},\tilde q_{\kern 1pt 6},
\tilde q_{\kern 1pt 7},\tilde q_{\kern 1pt 8},\tilde q_{\kern 1pt 9},
\tilde q_{\kern 1pt 10},\tilde q_{\kern 1pt 11},\tilde q_{\kern 1pt 12},
\tilde q_{\kern 1pt 13},\tilde q_{\kern 1pt 14}\bigr>.
\mytag{1.15}
$$
The polynomials $q_{\kern 1pt 1}$ through $q_{\kern 1pt 14}$ from \mythetag{1.15}
were calculated in \mycite{46} with the use of the Maxima symbolic computations 
package, but explicit formulas for them were not given. These formulas were 
given later in Appendix to \mycite{69}. For the theory of Gr\"obner bases and their
applications the reader is referred to \mycite{70}.\par
     Returning to the factor equations \mythetag{1.3} through \mythetag{1.10}
and noting that their left hand sides are multisymmetric polynomials, one can apply
Theorem~\mythetheorem{1.2} to them. This was done in \mycite{69} and the 
$E$-forms\footnotemark\ of the factor equations were derived in \mycite{69}. Here are 
the transformed factor equations:
\footnotetext{\ The E-form of a polynomial is its preimage under the mapping 
\mythetag{1.13}.}
$$
\allowdisplaybreaks
\gather
\hskip -2em
E_{10}^2-2\,E_{20}-L^2=0,\qquad
\mytag{1.16}\\
\vspace{2ex}
\hskip -2em
2\,E_{02}-4\,E_{20}-E_{01}^2+2\,E_{10}^2=0,\qquad
\mytag{1.17}\\
\vspace{2ex}
\hskip -2em
E_{10}\,E_{11}-3\,E_{03}-E_{21}+3\,E_{01}\,E_{02}-E_{20}\,E_{01}-E_{01}^3=0,\qquad
\mytag{1.18}\\
\vspace{2ex}
\hskip -2em
E_{01}\,E_{11}-E_{12}-3\,E_{30}+E_{10}\,E_{02}+E_{20}\,E_{10}-E_{01}^2\,E_{10}=0,
\qquad
\mytag{1.19}\\
\vspace{2ex}
\hskip -2em
\aligned
-E_{10}\,&E_{21}-E_{01}\,E_{12}-E_{01}\,E_{30}-E_{01}^3\,E_{10}+E_{01}^2\,E_{11}\,-\\
&-\,E_{02}\,E_{11}+E_{11}\,E_{20}-E_{10}\,E_{03}+2\,E_{10}\,E_{01}\,E_{02}=0.
\endaligned\qquad
\mytag{1.20}\\
\vspace{2ex}
\hskip -2em
\gathered
4\,E_{01}\,E_{10}\,E_{11}-3\,E_{01}^2\,E_{10}^2+2\,E_{10}^2\,E_{02}
+2\,E_{20}\,E_{01}^2-2\,E_{10}\,E_{12}\,-\\
-\,2\,E_{02}\,E_{20}-2\,E_{01}\,E_{21}-E_{11}^2-12\,E_{10}\,E_{30}+6\,E_{20}^2=0,
\endgathered\qquad
\mytag{1.21}\\
\vspace{2ex}
\hskip -2em
\gathered
4\,E_{01}\,E_{10}\,E_{11}-4\,E_{10}^2\,E_{02}-4\,E_{20}\,E_{01}^2-2\,E_{10}\,E_{12}
+10\,E_{02}\,E_{20}\,-\\
-\,2\,E_{01}\,E_{21}-E_{11}^2-12\,E_{01}\,E_{03}-3\,E_{01}^4-6\,E_{02}^2
+12\,E_{01}^2\,E_{02}=0,
\endgathered\qquad
\mytag{1.22}\\
\vspace{2ex}
\gathered
9\,E_{01}\,E_{03}\,E_{20}-7\,E_{01}^2\,E_{02}\,E_{20}+2\,E_{02}\,E_{10}\,E_{12}
-2\,E_{01}^2\,E_{10}\,E_{12}\,+\\
+\,3\,E_{03}\,E_{10}\,E_{11}+4\,E_{01}^3\,E_{10}\,E_{11}
-7\,E_{01}\,E_{02}\,E_{10}\,E_{11}-6\,E_{01}\,E_{03}\,E_{10}^2\,+\\
+\,8\,E_{01}^2\,E_{02}\,E_{10}^2+3\,E_{01}\,E_{11}\,E_{30}-2\,E_{01}\,E_{20}\,E_{21}
+E_{10}\,E_{12}\,E_{20}\,-\\
-\,E_{02}\,E_{10}^2\,E_{20}+E_{01}\,E_{10}\,E_{11}\,E_{20}+9\,E_{02}\,E_{10}\,E_{30}
-2\,E_{02}\,E_{20}^2\,+\\
+\,2\,E_{01}^2\,E_{20}^2-E_{11}^2\,E_{20}-3\,E_{12}\,E_{30}+E_{02}\,E_{11}^2
-E_{01}^2\,E_{11}^2\,-\\
-\,2\,E_{02}^2\,E_{10}^2+2\,E_{01}^4\,E_{20}+2\,E_{02}^2\,E_{20}-3\,E_{03}\,E_{21}\,-\\
-\,2\,E_{01}^3\,E_{21}+5\,E_{01}\,E_{02}\,E_{21}-6\,E_{01}^2\,E_{10}\,E_{30}
-3\,E_{01}^4\,E_{10}^2=0.
\endgathered\quad
\mytag{1.23}
\endgather
$$
In \mycite{69} the factor equations \mythetag{1.16}, \mythetag{1.17}, \mythetag{1.18}, 
\mythetag{1.19}, \mythetag{1.20}, \mythetag{1.21}, \mythetag{1.22}, and
\mythetag{1.23} were complemented with fourteen kernel equations 
$$
\xalignat 3
&\hskip -2em
\tilde q_{\kern 1pt 1}=0,
&&\tilde q_{\kern 1pt 2}=0,
&&\tilde q_{\kern 1pt 3}=0,\\
&\hskip -2em
\tilde q_{\kern 1pt 4}=0,
&&\tilde q_{\kern 1pt 5}=0,
&&\tilde q_{\kern 1pt 6}=0,\\
&\hskip -2em
\tilde q_{\kern 1pt 7}=0,
&&\tilde q_{\kern 1pt 8}=0,
&&\tilde q_{\kern 1pt 9}=0,
\mytag{1.24}\\
&\hskip -2em
\tilde q_{\kern 1pt 10}=0,
&&\tilde q_{\kern 1pt 11}=0,
&&\tilde q_{\kern 1pt 12}=0,\\
&\hskip -2em
\tilde q_{\kern 1pt 13}=0,
&&\tilde q_{\kern 1pt 14}=0.
\endxalignat
$$
As a result a huge system of twenty two polynomial equations with respect to ten
variables was obtained. In \mycite{69} this system was analyzed and luckily was 
reduced to a single equation. Here is this equation: 
$$
\hskip -2em
(2\,E_{11})^2+(E_{01}^2+L^2-E_{10}^2)^2-8\,E_{01}^2\,L^2=0. 
\mytag{1.25}
$$\par
     Equation \mythetag{1.25} turns out to be explicitly solvable, which is also 
a lucky event. Its general solution was discovered in \mycite{71}. This general 
solution of \mythetag{1.25}, upon eliminating $L$ by homogeneity, consists of one 
two-parameter solution and a series of one-parameter solutions. The main goal 
of the present paper is to propagate these explicit solutions back to the equations
\mythetag{1.16} through \mythetag{1.24} and formulate polynomial inverse problems
relating these solutions with perfect cuboids.\par
\head
2. Rational perfect cuboids. 
\endhead
     Let $x_1$, $x_2$, $x_3$, $d_1$, $d_2$, $d_3$ be edges and face diagonals of 
some perfect cuboid and let $L$ be its space diagonal. Then, dividing these numbers 
by $L$, we get a cuboid whose edges and face diagonals are given by rational numbers, 
while the space diagonal is equal to unity. Such a cuboid is called a rational 
perfect cuboid. Conversely, if we have a rational perfect cuboid with unit space
diagonal, we can take the common denominator of the rational numbers $x_1$, $x_2$, 
$x_3$, $d_1$, $d_2$, $d_3$ for $L$ and then, multiplying these numbers by $L$, we get
a perfect cuboid with integer edges and face diagonals whose space diagonal is equal 
to the integer number $L$. Thus, integer perfect cuboids and rational perfect cuboids 
with unit space diagonal are equivalent to each other.\par
     The equivalence of perfect cuboids and rational perfect cuboids was already 
used in \mycite{41} for deriving three cuboid conjectures. \pagebreak 
In the present paper we use this fact by setting $L=1$ and thus reducing the number 
of variables in the equations \mythetag{1.16} through \mythetag{1.23} and in the equation 
\mythetag{1.25}.\par
\head
3. The two-parameter case.
\endhead
     Let's substitute $L=1$ into \mythetag{1.25}. As a result, using the notation 
$E_{11}=x$, $E_{01}=y$, $E_{10}=z$, we get the equation coinciding with the equation 
\thetag{1.1} in \mycite{71}:
$$
\hskip -2em
(2\,E_{11})^2+(E_{01}^2+1-E_{10}^2)^2=8\,E_{01}^2. 
\mytag{3.1}
$$
This equation was solved in \mycite{71}. Theorem~2.1 from \mycite{71} yields the 
following two-parameter solution of the equation \mythetag{3.1}:
$$
\align
&\hskip -2em
E_{11}=-\frac{b\,(c^2+2-4\,c)}{b^2\,c^2+2\,b^2-3\,b^2\,c+c-b\,c^2\,+2\,b},
\mytag{3.2}\\
\vspace{1ex}
&\hskip -2em
E_{01}=-\frac{b(c^2+2-2\,c)}{b^2\,c^2+2\,b^2-3\,b^2\,c+c-b\,c^2+2\,b},
\mytag{3.3}\\
\vspace{1ex}
&\hskip -2em
E_{10}=-\frac{b^2\,c^2+2\,b^2-3\,b^2\,c\,-c}{b^2\,c^2+2\,b^2-3\,b^2\,c
+c-b\,c^2+2\,b}.
\mytag{3.4}
\endalign
$$
The numbers $b$ and $c$ are two rational parameters in \mythetag{3.2}, 
\mythetag{3.3}, and \mythetag{3.4}.\par
     Apart from $E_{11}$, $E_{10}$, and $E_{01}$ there are six other variables 
$E_{20}$, $E_{30}$, $E_{02}$, $E_{03}$, $E_{21}$, $E_{12}$ in the equations 
\mythetag{1.16} through \mythetag{1.23}. These variables are expressed through 
$E_{11}$, $E_{10}$, $E_{01}$, and $L$ by means of the formulas \thetag{4.1}, 
\thetag{4.3}, \thetag{4.6}, \thetag{4.7}, \thetag{5.1}, and \thetag{5.2} from 
\mycite{69}. Here are the formulas for $E_{20}$ and $E_{02}$:
$$
\xalignat 2
&\hskip -2em
E_{20}=\frac{1}{2}\,E_{10}^2-\frac{1}{2}\,L^2,
&&E_{02}=\frac{1}{2}\,E_{01}^2-L^2.
\mytag{3.5}
\endxalignat
$$
Substituting $L=1$, \mythetag{3.3}, and \mythetag{3.4} into \mythetag{3.5},
we derive
$$
\gather
\hskip -2em
\gathered
E_{20}=\frac{b}{2}\,(b\,c^2-2\,c-2\,b)\,(2\,b\,c^2-c^2-6\,b\,c+2
+4\,b)\,\times\\
\times\,(b\,c-1-b)^{-2}\,(b\,c-c-2\,b)^{-2},
\endgathered
\mytag{3.6}\\
\vspace{1ex}
\hskip -2em
\gathered
E_{02}=\frac{1}{2}\,(28\,b^2\,c^2-16\,b^2\,c-2\,c^2-4\,b^2-b^2\,c^4\,+\\
+\,4\,b^3\,c^4-12\,b^3\,c^3+4\,b\,c^3+24\,b^3\,c-8\,b\,c-2\,b^4\,c^4\,+\\
+\,12\,b^4\,c^3-26\,b^4\,c^2-8\,b^2\,c^3+24\,b^4\,c-16\,b^3-8\,b^4)\,\times\\
\times\,(b\,c-1-b)^{-2}\,(b\,c-c-2\,b)^{-2}.
\endgathered
\mytag{3.7}
\endgather
$$
The formulas for $E_{21}$ and $E_{12}$ are taken from \thetag{5.1} and \thetag{5.2} 
in \mycite{69} respectively:
$$
\allowdisplaybreaks
\align
&\hskip -2em
\aligned
E_{21}&=\frac{2\,E_{10}^3\,E_{11}+2\,E_{01}^2\,E_{10}\,E_{11}
-E_{01}\,E_{10}^4+E_{01}^5}
{8\,(E_{01}^2+E_{10}^2)\vphantom{\vrule height 11pt}}\,+\\
\vspace{1ex}
&+\,\frac{6\,E_{10}\,E_{11}\,L^2-2\,E_{01}\,E_{10}^2\,L^2
-8\,E_{01}^3\,L^2+3\,E_{01}\,L^4}
{8\,(E_{01}^2+E_{10}^2)\vphantom{\vrule height 11pt}},
\endaligned
\mytag{3.8}\\
\vspace{2ex}
&\hskip -2em
\aligned
E_{12}=&\frac{E_{01}^4\,E_{10}-2\,E_{01}^3\,E_{11}
-2\,E_{01}\,E_{10}^2\,E_{11}-E_{10}^5}
{8\,(E_{01}^2+E_{10}^2)\vphantom{\vrule height 11pt}}\,+\\
\vspace{1ex}
&\kern 5em +\,\frac{6\,E_{10}^3\,L^2-6\,E_{01}\,E_{11}\,L^2+3\,E_{10}\,L^4}
{8\,(E_{01}^2+E_{10}^2)\vphantom{\vrule height 11pt}}.
\endaligned
\mytag{3.9}
\endalign
$$
Substituting $L=1$, \mythetag{3.2}, \mythetag{3.3}, and \mythetag{3.4} into 
\mythetag{3.8} and \mythetag{3.9}, we derive
$$
\gather
\gathered
E_{21}=\frac{b}{2}\,(5\,c^6\,b-2\,c^6\,b^2+52\,c^5\,b^2-16\,c^5\,b
-2\,c^7\,b^2+2\,b^4\,c^8\,+\\
+\,142\,b^4\,c^6-26\,b^4\,c^7-426\,b^4\,c^5-61\,b^3\,c^6+100\,b^3\,c^5
+14\,c^7\,b^3\,-\\
-\,c^8\,b^3-20\,b\,c^2-8\,b^2\,c^2-16\,b^2\,c-128\,b^2\,c^4-200\,b^3\,c^3\,+\\
+\,244\,b^3\,c^2+32\,b\,c^3-112\,b^3\,c+768\,b^4\,c^4-852\,b^4\,c^3
+568\,b^4\,c^2\,+\\
+\,104\,b^2\,c^3-208\,b^4\,c+8\,c^4-4\,c^3+16\,b^3+32\,b^4-2\,c^5)\,\times\\
\times\,(b^2\,c^4-6\,b^2\,c^3+13\,b^2\,c^2-12\,b^2\,c+4\,b^2+c^2)^{-1}\,\times\\
\times\,(b\,c-1-b)^{-2}\,(b\,c-c-2\,b)^{-2},
\endgathered\qquad\quad
\mytag{3.10}\\
\vspace{1ex}
\hskip -2em
\gathered
E_{12}=(16\,b^6+32\,b^5-6\,c^5\,b^2+2\,c^5\,b-62\,b^5\,c^6+62\,b^6\,c^6\,-\\
-\,180\,b^6\,c^5+18\,b^5\,c^7-12\,b^6\,c^7-2\,b^5\,c^8+b^6\,c^8+248\,b^5\,c^2\,+\\
+\,248\,b^6\,c^2-96\,b^6\,c+321\,b^6\,c^4-180\,b^5\,c^3-144\,b^5\,c
-360\,b^6\,c^3\,+\\
+\,b^4\,c^8+8\,b^4\,c^6-6\,b^4\,c^7+18\,b^4\,c^5+7\,b^3\,c^6+90\,b^5\,c^5
-14\,b^3\,c^5\,-\\
-\,c^7\,b^3+17\,b^2\,c^4+28\,b^3\,c^3-28\,b^3\,c^2-4\,b\,c^3+8\,b^3\,c
-57\,b^4\,c^4\,+\\
+\,36\,b^4\,c^3+32\,b^4\,c^2-12\,b^2\,c^3-48\,b^4\,c-c^4+16\,b^4)\,\times\\
\times\,(b^2\,c^4-6\,b^2\,c^3+13\,b^2\,c^2-12\,b^2\,c+4\,b^2+c^2)^{-1}\,\times\\
\times\,(b\,c-1-b)^{-2}\,(b\,c-c-2\,b)^{-2}.
\endgathered\qquad
\mytag{3.11}
\endgather
$$
The formulas for $E_{30}$ and $E_{03}$ are taken from \thetag{4.6} and \thetag{4.7} 
in \mycite{69}:
$$
\align
&\hskip -2em
E_{30}=-\frac{1}{3}\,E_{12}-\frac{1}{6}\,E_{10}\,E_{01}^2
-\frac{1}{2}\,E_{10}\,L^2+\frac{1}{6}\,E_{10}^3
+\frac{1}{3}\,E_{01}\,E_{11},\quad
\mytag{3.12}\\
\vspace{1ex}
&\hskip -2em
E_{03}=-\frac{1}{3}\,E_{21}-\frac{1}{6}\,E_{01}\,E_{10}^2-\frac{5}{6}\,E_{01}\,L^2
+\frac{1}{6}\,E_{01}^3+\frac{1}{3}\,E_{10}\,E_{11}.\quad
\mytag{3.13}
\endalign
$$
The formulas \mythetag{3.12} and \mythetag{3.13} comprise $E_{21}$ and  $E_{12}$
from \mythetag{3.10} and \mythetag{3.11}. Applying $L=1$, \mythetag{3.10}, 
\mythetag{3.11}, \mythetag{3.2}, \mythetag{3.3}, and \mythetag{3.4} to 
\mythetag{3.12} and \mythetag{3.13}, we get
$$
\allowdisplaybreaks
\gather
\hskip -2em
\gathered
E_{03}=\frac{b}{2}\,(b^2\,c^4-5\,b^2\,c^3+10\,b^2\,c^2-10\,b^2\,c+4\,b^2+2\,b\,c\,+\\
+\,2\,c^2-b\,c^3)\,(2\,b^2\,c^4-12\,b^2\,c^3+26\,b^2\,c^2-24\,b^2\,c\,+\\
+\,8\,b^2-c^4\,b+3\,b\,c^3-6\,b\,c+4\,b+c^3-2\,c^2+2\,c)\,\times\\
\times\,((b^2\,c^4-6\,b^2\,c^3+13\,b^2\,c^2-12\,b^2\,c+4\,b^2+c^2)^{-1}\,\times\\
\times\,(b\,c-1-b)^{-2}\,(-c+b\,c-2\,b)^{-2},
\endgathered
\mytag{3.14}\\
\vspace{1ex}
\hskip -2em
\gathered
E_{30}=c\,b^2\,(1-c)\,(c-2)\,(b\,c^2-4\,b\,c+2+4\,b)\,\times\\
\times\,(2\,b\,c^2-c^2-4\,b\,c+2\,b)\,\times\\
\times\,(b^2\,c^4-6\,b^2\,c^3+13\,b^2\,c^2-12\,b^2\,c+4\,b^2+c^2)^{-1}\,\times\\
\times\,(b\,c-1-b)^{-2}\,(-c+b\,c-2\,b)^{-2}.
\endgathered
\mytag{3.15}
\endgather
$$\par
     Note that the quantities $E_{10}$, $E_{20}$, and $E_{30}$ given by the formulas 
\mythetag{3.4}, \mythetag{3.6},	\mythetag{3.15} are the values of three elementary 
multisymmetric polynomials in the left column of \mythetag{1.11}. One can easily see 
that these polynomials coincide with the regular symmetric polynomials of three 
variables $x_1$, $x_2$, and $x_3$ (see \mycite{72}). Therefore, if we write the cubic 
equation  $(x-x_1)(x-x_2)(x-x_3)=0$, it expands to
$$
\hskip -2em
x^3-E_{10}\,x^2+E_{20}\,x-E_{30}=0.
\mytag{3.16}
$$
Similarly, the quantities $E_{01}$, $E_{02}$, and $E_{03}$ given by the formulas 
\mythetag{3.3}, \mythetag{3.7},	and \mythetag{3.14} are the values of the elementary 
multisymmetric polynomials in the right column of \mythetag{1.11}. These three 
polynomials coincide with regular symmetric polynomials of three variables $d_1$, 
$d_2$, and $d_3$. For this reason, if we write the cubic equation 
$(d-d_1)(d-d_2)(d-d_3)=0$, this equation expands to
$$
\hskip -2em
d^{\kern 1pt 3}-E_{01}\,d^{\kern 1pt 2}+E_{02}\,d-E_{03}=0.
\mytag{3.17}
$$
The remaining three quantities $E_{11}$, $E_{21}$, and $E_{12}$ given by the formulas 
\mythetag{3.2}, \mythetag{3.10}, \mythetag{3.11} are the values of three elementary
multisymmetric polynomials in \mythetag{1.12}. They lead to the following auxiliary
polynomial equations:
$$
\hskip -2em
\aligned
&x_1\,x_2\,d_3+x_2\,x_3\,d_1+x_3\,x_1\,d_2=E_{21},\\
&x_1\,d_2+d_1\,x_2+x_2\,d_3+d_2\,x_3+x_3\,d_1+d_3\,x_1=E_{11},\\
&x_1\,d_2\,d_3+x_2\,d_3\,d_1+x_3\,d_1\,d_2=E_{12}.
\endaligned
\mytag{3.18}
$$
Now we can formulate the following two inverse cuboid problems.
\myproblem{3.1} Find all pairs of rational numbers $b$ and $c$ for which the
cubic equations \mythetag{3.16} and \mythetag{3.17} with the coefficients given
by the formulas \mythetag{3.4}, \mythetag{3.6},	\mythetag{3.15}, \mythetag{3.3}, 
\mythetag{3.7},	\mythetag{3.14} possess positive rational roots $x_1$, $x_2$, 
$x_3$, $d_1$, $d_2$, $d_3$ obeying the auxiliary polynomial equations 
\mythetag{3.18} whose right hand sides are given by the formulas \mythetag{3.2}, 
\mythetag{3.10}, \mythetag{3.11}. 
\endproclaim
\myproblem{3.2} Find at least one pair of rational numbers $b$ and $c$ for which 
the cubic equations \mythetag{3.16} and \mythetag{3.17} with the coefficients given
by the formulas \mythetag{3.4}, \mythetag{3.6},	\mythetag{3.15}, \mythetag{3.3}, 
\mythetag{3.7},	\mythetag{3.14} possess positive rational roots $x_1$, $x_2$, 
$x_3$, $d_1$, $d_2$, $d_3$ obeying the auxiliary polynomial equations 
\mythetag{3.18} whose right hand sides are given by the formulas \mythetag{3.2}, 
\mythetag{3.10}, \mythetag{3.11}. 
\endproclaim
\head
4. The first one-parameter case. 
\endhead
     The four one-parameter cases considered below correspond to the special 
solutions of Theorem 2.1 in \mycite{71}. The first one-parameter case 
corresponds to the following choice of sign options in this theorem:
$$
\xalignat 2
&\hskip -2em
x=y,
&&y=z+1.
\mytag{4.1}
\endxalignat
$$
Equations \mythetag{4.1} are easily resolved in the 
one-parameter form:
$$
\xalignat 3
&\hskip -2em
x=c,
&&y=c,
&&z=c-1.
\mytag{4.2}
\endxalignat
$$
Taking into account the notations $E_{11}=x$, $E_{01}=y$, $E_{10}=z$, we derive
the following one-parameter solution of the equation \mythetag{3.1} from
\mythetag{4.2}:	
$$
\xalignat 3
&\hskip -2em
E_{11}=c,
&&E_{01}=c,
&&E_{10}=c-1.
\mytag{4.3}
\endxalignat
$$
Substituting $L=1$ and \mythetag{4.3} into \mythetag{3.5}, \mythetag{3.8}, 
\mythetag{3.9}, \mythetag{3.12}, and \mythetag{3.13}, we get
$$
\xalignat 2
&\hskip -2em
E_{20}=\frac{c^2-2\,c}{2},
&&E_{02}=\frac{c^2-2}{2},\\
&\hskip -2em
E_{21}=\frac{c^2-2\,c}{2},
&&E_{12}=1,
\mytag{4.4}\\
&\hskip -2em
E_{30}=0,
&&E_{03}=\frac{c^2-2\,c}{2}.
\endxalignat
$$
In the formulas \mythetag{4.4} we find $E_{30}=0$, i\.\,e\. the last term $E_{30}$
of the cubic equation \mythetag{3.16} is zero. Then one of its roots $x_1$, $x_2$, 
or $x_3$ is zero. But a non-degenerate cuboid cannot have a zero edge. As a result 
we have the following theorem.
\mytheorem{4.1} There are no perfect cuboids associated with the one-parameter 
solution \mythetag{4.3} of the equation \mythetag{3.1}. 
\endproclaim
\head
5. The second one-parameter case. 
\endhead
    In the second one-parameter case we choose the following sign options
in Theorem 2.1 from \mycite{71}: $x=y$ and $y=-z-1$. Then instead of \mythetag{4.3} 
we derive 
$$
\xalignat 3
&\hskip -2em
E_{11}=c,
&&E_{01}=c,
&&E_{10}=-c-1.
\mytag{5.1}
\endxalignat
$$
Substituting $L=1$ and \mythetag{5.1} into \mythetag{3.5}, \mythetag{3.8}, 
\mythetag{3.9}, \mythetag{3.12}, and \mythetag{3.13}, we get
$$
\xalignat 2
&\hskip -2em
E_{20}=\frac{c^2+2\,c}{2},
&&E_{02}=\frac{c^2-2}{2},\\
&\hskip -2em
E_{21}=-\frac{c^2+2\,c}{2},
&&E_{12}=1,
\mytag{5.2}\\
&\hskip -2em
E_{30}=0,
&&E_{03}=-\frac{c^2+2\,c}{2}.
\endxalignat
$$
In \mythetag{5.2} we again see $E_{30}=0$. Therefore we can formulate the
following theorem. 
\mytheorem{5.1} There are no perfect cuboids associated with the one-parameter 
solution \mythetag{5.1} of the equation \mythetag{3.1}.
\endproclaim
\head
6. The third and the fourth one-parameter cases. 
\endhead
     These two one-parameter cases correspond to the last two sign options 
in Theorem 2.1 from \mycite{71}. \pagebreak 
These two sign options are given by the formulas $x=-y$ and $y=\pm\,(z+1)$. 
They yield the following one-parameter solutions for \mythetag{3.1}:
$$
\xalignat 3
&\hskip -2em
E_{11}=c,
&&E_{01}=-c,
&&E_{10}=\pm\,c-1.
\mytag{6.1}
\endxalignat
$$
The formulas $E_{11}=c$ and $E_{01}=-c$ mean that $E_{11}$ and $E_{01}$ are of
opposite signs or both of them are zero. On the other hand, we have two equalities
$$
\align
&x_1\,d_2+d_1\,x_2+x_2\,d_3+d_2\,x_3+x_3\,d_1+d_3\,x_1=E_{11},\\
&d_1+d_2+d_3=E_{01}
\endalign
$$
derived from \mythetag{1.11} and \mythetag{1.12}. The left hand sides of these
equalities are positive numbers since edges and face diagonals of a perfect cuboid 
are positive. The contradiction obtained leads to the following theorem.
\mytheorem{6.1} There are no perfect cuboids associated with the one-parameter 
solutions \mythetag{6.1} of the equation \mythetag{3.1}. 
\endproclaim
\head
7. Conclusions. 
\endhead
    Theorems~\mythetheorem{4.1}, \mythetheorem{5.1}, and \mythetheorem{6.1}
show that the inverse problem~\mytheproblem{3.1} formulated in Section 3 is 
equivalent to finding all perfect cuboids, while the problem~\mytheproblem{3.2}
is equivalent to finding at least one perfect cuboid. 

\enddocument
\end